\documentclass{article}

\usepackage{graphicx}
\usepackage{epsfig}
\usepackage{amssymb}
\usepackage{amsmath}
\usepackage{amsfonts, latexsym,cite}
\usepackage{stmaryrd}

\newcommand{\ds}{\displaystyle}
\newcommand{\be}{\begin{equation}}

\newcommand{\ee}{\end{equation}}
\newcommand{\ba}{\begin{array}}
\newcommand{\ea}{\end{array}}
\newcommand{\bea}{\begin{eqnarray}}
\newcommand{\eea}{\end{eqnarray}}
\newcommand{\bi}{\begin{itemize}}
\newcommand{\ei}{\end{itemize}}

\newtheorem{prop}{Proposition}

\usepackage{authblk}

\begin{document}

\title{A higher rank extension of the Askey-Wilson Algebra }
\author{S. Post and A. Walter}
\affil[]{University of Hawai`i at M\-{a}noa\\
Department of Mathematics\\
2565 McCarthy Mall\\
Honolulu, HI 96815 }
\maketitle
\abstract{A novel generalization of the Askey-Wilson algebra is presented and shown to be associated with coproducts in the quantum algebra $U_q(su(1,1))$. This algebra has 15 non-commuting generators given by $Q^{(A)}$, with $A\subset \{1,2,3,4\}$ and their 5 linearly independent inversions generated by the algebra automorphism $q\rightarrow q^{-1}$, $E\leftrightarrow F.$ The set of generators  can be split into operators fixed under inversion, and those with an orientation under this inversion. We then show that the generators will either commute or satisfy q-commutator relations linear in the generators, with the restriction that parity operators commute with generators of opposite parity and the q-commutation relations are between those with the same parity. Finally, we give a novel algebra expression satisfied by the generators involving only the natural generators, i.e. those arising from the coupling scheme. 
 
	  }

\section{Introduction}
The connection between orthogonal polynomials and symmetry algebras has long been a source of  fruitful research. For example, many orthogonal polynomials and their properties can be described in terms of Lie algebras, see for example W. Miller's monograph, Lie Theory and Special Functions \cite{Miller68}. However, many special functions and orthogonal polynomials were not completely characterized or obtained in all generality within the framework of Lie algebras and Lie groups. 

In 1992\cite{Zhedanov1991}, A. Zhedanov extended this approach to quadratic algebras showing that the Askey-Wilson polynomials, the most generic bispectral hypergeometic polynomials, give rise to a quadratic algebra, dubbed $AW(3).$ This algebra is generated by the action of operators associated with multiplication by the polynomial's variable and the difference operator for which the polynomials are eigenfunctions. Later, this algebraic structure was generalized to different classes in the Askey Scheme (Racah, Hahn, etc. )\cite{granovskii1992mutual}. The fact that these polynomials are dual in the sense of Leonard's theorem arises from their representation as interbasis expansion between eigenvectors of a Leonard pair. In the study of these systems, P. Terwilliger, including work with Vidunas,  was able to show that any such Leonard pair gives rise to exactly these algebra relations and discuss their connection with the polynomial systems \cite{Terwilliger2001, TerwilligerVidunas2004, Terwilliger2006}. The study of these quadratic algebras has found many fruitful applications including to the superintegrable systems\cite{KKM2005, Post11, KMP2013, genest2014superintegrability}.   

From a slightly different perspective, the Askey-Wilson polynomials and their degeneration were being studied in the context of quantum groups and in particular as 6j or Racah expansions of the quantum algebras $U_{q}(\mathfrak{su}(2))$ and $U_q(\mathfrak{su}(1,1))$ \cite{kirillow1989representations, koelink1998convolutions}. This was followed by extensions into higher dimensions by considering 3nj expansions \cite{van20033nj, groenevelt2011quantum }. Recently,   multivariate generalizations of the q-Racah polynomials were shown explicitly to arise as 3nj coefficients for the algebra    $U_q(\mathfrak{su}(1,1))$ \cite{genest2017coupling}.

In this work, we use the explicit connection between the two approaches to construct a higher-rank generalization of the Askey-Wilson algebra. In the triple product case, it has been shown that the Racah problem for $sl_q(2)$, supports a representation of the algebra $AW(3)$  \cite{GranZhed1993}. This problem has more recently been studies in depth including deriving decomposition rules \cite{Huang2017} and the connection with matrix models and conformal field theory \cite{MMS2017}. In this paper, we extend the problem to the quadruple product to obtain a novel generalization of the Askey-Wilson algebra, which we call $AW(4).$  It should be noted that this connection has already been established for coupling schemes of $su(1,1)$, bivariate Racah polynomials, and higher rank generalizations of the quadratic Racah algebra\cite{genest2013equitable, post2015racah, de1610higher}.

\section{$AW(3)$ algebra and the algebra $U_q(su(1,1))$}
In this section, we will describe the well-known Askey-Wilson algebra in terms of  decomposition schemes of tensor products of irreducible representations of the quantum algebra $U_q(su(1,1)).$ 
\subsection{The algebra $U_q(su(1,1))$ and its coproduct structure}
Given $q\in \mathbb{C}\setminus \{0, \pm 1\}$  the algebra $U_q(sl(2,\mathbb{C}))$ is the associative, complex algebra with generators $E,F,K,K^{-1}$ subject to the relations
\be KK^{-1}=1=K^{-1}K,\qquad KE=qEK, \qquad qKF=FK, \qquad [E,F]=\frac{K^2-K^{-2}}{q-q^{-1}}. \ee
Here the Casimir element
\be \displaystyle \Omega=\frac{q^{-1}K^2+qK^{-2}-2}{(q-q^{-1})^2}+EF \ee
is a central element, and generates the center if $q$ is not a root of unity. 
There is a comultiplication $\Delta:U_q(sl(2,\mathbb{C}))\to U_q(sl(2,\mathbb{C}))\otimes U_q(sl(2,\mathbb{C}))$ defined on the generators by
\be \Delta(K)=K\otimes K,\qquad \Delta(K^{-1})=K^{-1}\otimes K^{-1}, \qquad \Delta(E)=K\otimes E + E \otimes K^{-1}, \qquad \Delta(F)=K\otimes F + F\otimes K^{-1} \ee
Note that this comultiplication is coassociative but not cocommutative.

If $q\in \mathbb{R}\setminus \{0,\pm 1\}$, we have a $*$-structure on $U_q(sl(2,\mathbb{C}))$ defined by
\be K^*=K,\qquad (K^{-1})^*=K^{-1}, \qquad E^*=-F, \qquad F^*=-E \ee
We denote the corresponding $*$-algebra by $U_q(su(1,1))$.

\subsection{Decomposition of tensor products}
Suppose we have a pair of positive discrete series representations with basis vectors $e_{n_i}^{k_i}$ satisfying 
\be \label{posdiscrete} \begin{gathered} \pi_{k_{i}}(K)e_{n_i}^{(k_i)}=q^{(k_i+n_i)}e_{n_i}^{(k_i)}, \qquad  \pi_{k_{i}}(K^{-1})e_{n_i}^{(k_i)}=q^{-(k_i+n_i)}e_{n_i}^{(k_i)}, \\ \pi_{k_{i}}(E)e_{n_i}^{(k_i)} = \frac{q^{-1/2 -k_i-n_i}}{q^{-1}-q}\sqrt{(1-q^{2n_i+2})(1-q^{4k_i+2n_i})}e_{n_i+1}^{(k_i)},\\ \pi_{k_{i}}(F)e_{n_i}^{(k_i)} = -\frac{q^{1/2 -k_i-n_i}}{q^{-1}-q}\sqrt{(1-q^{2n_i})(1-q^{4k_i+2n_i-2})}e_{n_i-1}^{(k_i)}.
\end{gathered} \ee
The tensor product of these representations supports a reducible representation of the algebra consistent with the co-product structure: 
\be \label{tensorproduct} \begin{gathered}
\pi_{k_{12}}(K)\equiv \pi_{k_{1}}(K) \otimes \pi_{k_2}(K), \qquad \pi_{k_{12}}(K^{-1})\equiv \pi_{k_{1}}(K^{-1}) \otimes \pi_{k_2}(K^{-1}),\\ 
\pi_{k_{12}}(E)\equiv \pi_{k_{1}}(E)\otimes \pi_{k_2}(K^{-1}) + \pi_{k_{1}}(K) \otimes \pi_{k_2}(E), \\ 
\pi_{k_{12}}(F)\equiv \pi_{k_{1}}(F)\otimes \pi_{k_2}(K^{-1}) + \pi_{k_{1}}(K) \otimes \pi_{k_2}(F).  \end{gathered} \ee
The action of the operators on the basis $e_{n_1}^{(k_1)} \otimes e_{n_2}^{(k_2)}$ is
\be \begin{gathered} \pi_{k_{12}}(K) e_{n_1}^{(k_1)} \otimes e_{n_2}^{(k_2)}= q^{(k_1+k_2+n_1+n_2)} e_{n_1}^{(k_1)} \otimes e_{n_2}^{(k_2)},\\ \pi_{k_{12}}(K^{-1}) e_{n_1}^{(k_1)} \otimes e_{n_2}^{(k_2)}= q^{-(k_1+k_2+n_1+n_2)} e_{n_1}^{(k_1)} \otimes e_{n_2}^{(k_2)}\\
\begin{split}
\pi_{k_{12}}(E) e_{n_1}^{(k_1)} \otimes e_{n_2}^{(k_2)}&=\frac{q^{-1/2 -k_1-n_1-k_2-n_2}}{q^{-1}-q} \sqrt{(1-q^{2n_1+2})(1-q^{4k_1+2n_1})} e_{n_1+1}^{(k_1)}\otimes e_{n_2}^{(k_2)}\\
&+ \frac{q^{-1/2 +k_1+n_1-k_2-n_2}}{q^{-1}-q} \sqrt{(1-q^{2n_2+2})(1-q^{4k_2+2n_2})}e_{n_1}^{(k_1)}\otimes e_{n_2+1}^{(k_2)},
\end{split}\\
\begin{split}
\pi_{k_{12}}(F) e_{n_1}^{(k_1)} \otimes e_{n_2}^{(k_2)}&=-\frac{q^{1/2 -k_1-n_1-k_2-n_2}}{q^{-1}-q} \sqrt{(1-q^{2n_1})(1-q^{4k_1+2n_1-2})} e_{n_1-1}^{(k_1)}\otimes e_{n_2}^{(k_2)} \\
& - \frac{q^{1/2 +k_1+n_1-k_2-n_2}}{q^{-1}-q} \sqrt{(1-q^{2n_2})(1-q^{4k_2+2n_2-2})}e_{n_1}^{(k_1)}\otimes e_{n_2-1}^{(k_2)}. 
\end{split}
\end{gathered} \ee
The product of two positive discrete series can be be decomposed into irreducible positive discrete series indexed by a positive integer $x$
\[ \pi_{k_1}\otimes \pi_{k_2}=\oplus_{x=0}^{\infty} \pi_{k_1+k_2 +x}.\]
The action of the Casimir operator 
\be \pi_{k_{12}}(\Omega) \equiv \frac{q^{-1}\pi_{k_{12}}(K)^2+q\pi_{k_{12}}(K^{-1})^2-2}{(q^{-1}-q)^2} +\pi_{k_{12}}(E)\pi_{k_{12}}(F),\ee
on any element $e_{n}^{(k_{12}+x)}$ in the irreducible component is
\be \pi_{k_{12}}(\Omega)e_{n}^{(k_{12}+x)}= \frac{q^{2k_{12}+2x-1}+q^{1-2k_{12}-2x}-2}{(q^{-1}-q)^2}e_{n}^{(k_{12}+x)}.\ee
Here we have use the notation $k_{12}=k_1+k_1, $ and later $k_{123}=k_1+k_2+k_3$. 

This process can be continued by taking the tensor product of these irreducible components and tensoring with an additional discrete series $e_{n_3}^{(k_3)}$ to obtain the following decomposition 
\[ \pi_{k_1}\otimes \pi_{k_2}\otimes \pi_{k_3}=\oplus_{0\leq x\leq N}\, \,  \pi_{k_{12} +x, k_{123} +N},\]
with basis vectors written $e_{n}^{(k_{12} +x, k_{123} +N)}.$
Of course the Casimir operator 
\[ \pi_{123}(\Omega) \equiv \frac{q^{-1}\pi_{k_{123}}(K)^2+q\, \pi_{k_{123}}(K^{-1})^2-2}{(q^{-1}-q)^2} +\pi_{k_{123}}(E)\pi_{k_{123}}(F),\] 
will be central in the algebra and so will act as a constant on any irreducible representation,
\be \pi_{k_{123}}(\Omega )e_{n}^{(k_{12}+x, k_{123}+N)}= \frac{q^{2k_{123}+2x-1}+q^{1-2k_{123}-2x}-2}{(q^{-1}-q)^2}e_{n}^{(k_{12}+x, k_{123}+N)}.\ee
Additionally, if we define an intermediate Casimir operator 
\[ Q^{(12)}\equiv \pi_{12}(\Omega)\otimes I,\]
this operators will also commute with the  generators of the $U_q(su(1,1))$ algebra
\be \label{pi123 def} \begin{gathered}
\pi_{k_{123}}(K)= \pi_{k_{12}}(K) \otimes \pi_{k_3}(K), \qquad \pi_{k_{123}}(K^{-1})= \pi_{k_{12}}(K^{-1}) \otimes \pi_{k_3}(K^{-1}),\\ 
\pi_{k_{123}}(E)= \pi_{k_{12}}(E)\otimes \pi_{k_3}(K^{-1}) + \pi_{k_{12}}(K) \otimes \pi_{k_3}(E), \\ 
\pi_{k_{123}}(F)= \pi_{k_{12}}(F)\otimes \pi_{k_3}(K^{-1}) + \pi_{k_{12}}(K) \otimes \pi_{k_3}(F).  \end{gathered} \ee

Notice that in constructing this triple tensor product we coupled the third component on the left. If we had instead coupled the component on the right, we would have obtained a different representation all together, as in
 \[ \pi_{k_{312}}(E)=\pi_{k_{3}}(K) \otimes \pi_{k_{12}}(E)+\pi_{k_{3}}(E)\otimes \pi_{k_{12}}(K^{-1}). \]
This is due to the fact that the coproduct structure is not cocommutative. It is however associative and so we can obtain the same tensor product by first coupling the representations $\pi_{k_2}$ and $\pi_{k_3}$ and then coupling to that the representation $\pi_{k_1}$ on the right. The operators then become 
\be  \begin{gathered}
\pi_{k_{123}}(K)= \pi_{k_{1}}(K) \otimes \pi_{k_{23}}(K), \qquad \pi_{k_{123}}(K^{-1})= \pi_{k_{1}}(K^{-1}) \otimes \pi_{k_{23}}(K^{-1}),\\ 
\pi_{k_{123}}(E)= \pi_{k_{1}}(E)\otimes \pi_{k_{23}}(K^{-1}) + \pi_{k_{1}}(K) \otimes \pi_{k_{23}}(E), \\ 
\pi_{k_{123}}(F)= \pi_{k_{1}}(F)\otimes \pi_{k_{23}}(K^{-1}) + \pi_{k_{1}}(K) \otimes \pi_{k_{23}}(F).  \end{gathered} \ee
It is straightforward to check that these operators agree with those given by (\ref{pi123 def}). 
From this coupling scheme, we obtain the operator
\[ Q^{(23)}= I \otimes \pi_{k_{23}}(\Omega),\]
which also commutes with the algebra generators and acts as a constant on the decomposition obtained by first decomposing the $\pi_{k_{23}}$ representation, i.e.
\[ \pi_{k_1}\otimes \pi_{k_2}\otimes \pi_{k_3}=\oplus_{0\leq y\leq N}\, \,  \pi_{k_{23} +y, k_{123} +N}.\]
The operators $Q^{(12)}$ and $Q^{(23)}$ will not however themselves commute and they generate the Askey-Wilson algebra $AW(3).$

\subsection{The algebra $AW(3)$}
Although the tensor products of irreducible representations is a useful motivation, it is not strictly necessary to define the algebra. The algebra can be defined directly from the following operators 
\[ Q^{(1)}=\Omega \otimes I \otimes I,\quad  Q^{(2)}=I \otimes\Omega  \otimes I,\quad  Q^{(3)}=I \otimes I \otimes \Omega,\]
\[ Q^{(12)}=\Delta \Omega \otimes I ,\quad  Q^{(23)}=I \otimes\Delta\Omega, \]
and 
\[ Q^{(123)}=\Delta(\Delta(\Omega)).\]
The Askey-Wilson algebra is the algebra generated by these operators which satisfy relations that can best be expressed via the q-commutator 
\be [A, B]_q\equiv q A B-q^{-1} B A.\ee
The algebra relations are
\be \label{QAW3}\begin{aligned}	[]	
[[Q^{(12)},Q^{(23)}]_q,Q^{(12)}]_q &=-2(Q^{(12)})^2-2\{Q^{(12)},Q^{(23)}\}+BQ^{(12)}+Q^{(23)}+D_1\\
[[Q^{(23)},Q^{(12)}]_q,Q^{(23)}]_q &=-2(Q^{(23)})^2-2\{Q^{(12)},Q^{(23)}\}+BQ^{(23)}+Q^{(12)}+D_2,\\
\end{aligned} \ee
where $B,D_1,D_2$ composed of central elements,
\begin{gather*}
B=(q-q^{-1})^2(Q^{(1)}Q^{(3)}+Q^{(2)}Q^{(123)})+2(Q^{(1)}+Q^{(2)}+Q^{(3)}+Q^{(123)})\\
D_1=2(Q^{(1)}Q^{(3)}+Q^{(2)}Q^{(123)}) -\frac{2q(Q^{(1)}+Q^{(2)}+Q^{(3)}+Q^{(123)})}{(q+1)^2}-(q+q^{-1})(Q^{(1)}Q^{(123)}+Q^{(2)}Q^{(3)})
+\frac{2q^2}{(q+1)^4},\\
D_2=2(Q^{(1)}Q^{(3)}+Q^{(2)}Q^{(123)}) -\frac{2q(Q^{(1)}+Q^{(2)}+Q^{(3)}+Q^{(123)})}{(q+1)^2}-(q+q^{-1})(Q^{(3)}Q^{(123)}+Q^{(1)}Q^{(2)})
+\frac{2q^2}{(q+1)^4}.\\
\end{gather*}
These algebra relations are equivalent to the algebra $AW(3)$ given in \cite{Zhedanov1991}, also referred to as the Leonard algebra\cite{TerwilligerUnpublished}. 
If we shift and scale the operators $\ds Q^{(\cdot)}\rightarrow - \frac{(q+q^{-1})Q^{(\cdot)}+2}{(q-q^{-1})^2}$ the relations become linear in the generators
\be \label{AW3} 
\begin{aligned}	[]	
[[Q^{(12)},Q^{(23)}]_q,Q^{(12)}]_q &=(q-q^{-1})^2(BQ^{(12)}+Q^{(23)}+Q^{(1)}Q^{(123)}+Q^{(2)}Q^{(3)})\\
[[Q^{(23)},Q^{(12)}]_q,Q^{(23)}]_q &=(q-q^{-1})^2(BQ^{(23)}+Q^{(12)}+Q^{(3)}Q^{(123)}+Q^{(1)}Q^{(2)})\\
&B=(Q^{(1)}Q^{(3)}+Q^{(2)}Q^{(123)})\\
\end{aligned}
\ee
We can then define another generator $Q^{(13)}$,
\[ Q^{(13)} \equiv \frac{1}{q-q^{-1}} [Q^{(12)}, Q^{(23)}]_q -Q^{(1)}Q^{(3)} -Q^{(2)}Q^{(123)}\]
Then this relationship can be expressed symmetrically as:
\be
\begin{aligned}	[]	
	\frac{1}{q-q^{-1}}[Q^{(12)},Q^{(23)}]_q &=Q^{(13)}+Q^{(3)}Q^{(1)}+Q^{(2)}Q^{(123)}\\
	\frac{1}{q-q^{-1}}[Q^{(13)},Q^{(12)}]_q &=Q^{(23)}+Q^{(2)}Q^{(3)}+Q^{(1)}Q^{(123)}\\
	\frac{1}{q-q^{-1}}[Q^{(23)},Q^{(13)}]_q &=Q^{(12)}+Q^{(1)}Q^{(2)}+Q^{(3)}Q^{(123)}\\
\end{aligned}
\ee

Note that it is perhaps better to refer to this presentation as the 'universal' Askey-Wilson algebra\cite{Terwilliger2011} as the operators $Q^{(j)}$ and $Q^{(123)}$ should be considered as central elements rather that complex constants. For the remainder of the paper, we consider only the shifted form of the operators. 

\section{An extension to higher dimensions: $AW(4)$}
\subsection{Decomposition Schemes}Consider now the extension to the tensor product of 4 irreducible representations. These can be coupled iteratively from left to right, generating the sequence of operators 
\[  Q^{(12)}, Q^{(123)}, Q^{(1234)}\]
 from right to left 
\[  Q^{(34)}, Q^{(234)}, Q^{(1234)},\]
or starting from the middle 
\[ Q^{(23)}, Q^{(123)}, Q^{(1234)},\qquad Q^{(23)}, Q^{(234)}, Q^{(1234)}.\]
Thus, we will have all operators of the form $Q^{(i)}$, $Q^{(ij)}$,  $Q^{(ijk)}$ and $Q^{(1234)}$, with each index consecutive and in increasing order. 
 This is in distinction to the classical case\cite{genest2013equitable, genest2014superintegrability, post2015racah, de1610higher}, in that the order of the indices is important and must be preserved to guarantee that the elements commute with the total Casimir.   

A decomposition scheme is then given by a set of 2 mutually commuting operators, in addition to the single Casimirs $Q^{(i)}$ and the total Casimir $Q^{(1234)}$ which form the center of the algebra.  There are four Jucys-Murphy type couplings as given above, and also a single pair of cylindrical type or 9-j type given by $Q^{(12)}$ and $Q^{(34)}.$ In\cite{genest2017coupling}, the interbasis expansions associated with the first and the last of the above sequence of operators was shown to give bivariate q-Racah polynomaials. 

\subsection{The algebra $AW(4)$}
We define this algebra, $AW(4),$ to be the algebra generated by these intermediate Casimir operators defined as 
\[ Q^{(i)} =\underbrace{I\otimes \ldots I }_{i-1}\otimes \Omega\otimes \underbrace{ I\otimes \ldots I}_{n-i},\]
\[ Q^{(ij)}=\underbrace{I\otimes \ldots I}_{i-1} \otimes \Delta(\Omega) \otimes\underbrace{I\otimes \ldots I}_{n-i-1},\]
\[ Q^{(ijk)}=\underbrace{I\otimes \ldots I}_{i-1}\otimes \Delta^2(\Omega)\otimes \underbrace{  I\otimes \ldots I}_{n-i-2},\]
\[ Q^{(1234)}=\Delta^3(\Omega). \]
As a reminder, we  use the shifted and scaled Casimir $\Omega=-(q^{-1}K^2+qK^{-2}+(q-q^{-1})^{2}EF)(q+q^{-1})^{-1}$ for a more compact presentation of the relations, as in the rank 3 case.
\begin{prop}
For ordered subsets $A,B \subset \{1, 2, 3, 4\}$, the operators $Q^{(A)}$ and $Q^{(B)}$ commute whenever $A\cap B=\emptyset$ or $A\subset B$. 
\end{prop}
The proof in the case of an empty intersection is straightforward from the definition of the operators, since each will act on disjoint components of the tensor product. If $A\subset B$ then the operator $Q^{(A)}$ will compute with the generators of the algebra $E^{(B)}$, $F^{(B)}$, $K^{(B)}$ and $(K^{(B)})^{-1}$. Thus it will commute with the Casimir operator $Q^{(B)}.$ 

The algebra AW(4) contains five copies of AW(3). Each corresponding to a triangle in Figure \ref{AW4figure} consisting of two solid edges and one of the dashed edges. Similar to the case of $AW(3)$ above,  we define 5 new generators from q-commutators of the original:
\be 
\begin{aligned} []
Q^{(13)} &\equiv \frac{1}{q-q^{-1}} [Q^{(12)}, Q^{(23)}]_q -Q^{(1)}Q^{(3)} -Q^{(2)}Q^{(123)}\\
Q^{(24)} &\equiv \frac{1}{q-q^{-1}} [Q^{(23)}, Q^{(34)}]_q -Q^{(2)} Q^{(4)} -Q^{(3)} Q^{(234)}\\
Q^{(124)} &\equiv \frac{1}{q-q^{-1}} [Q^{(34)}, Q^{(123)}]_q -Q^{(12)} Q^{(4)} -Q^{(3)} Q^{(1234)}\\
Q^{(14)} &\equiv \frac{1}{q-q^{-1}}[Q^{(123)}, Q^{(234)}]_q -Q^{(1)}Q^{(4)} -Q^{(23)}Q^{(1234)}\\
Q^{(134)} &\equiv \frac{1}{q-q^{-1}} [Q^{(234)}, Q^{(12)}]_q -Q^{(1)} Q^{(34)} -Q^{(2)}Q^{(1234)}\\
\end{aligned}
\ee
Next we  introduce the involution $I$ which sends $q\to q^{-1}$ and $E\leftrightarrow F$. This will give us $IQ^{(A)}=Q^{(A )}$ for all of the original generators, and for $A=13,134,14,124,24$, $IQ^{(A)}$ is the element which looks like $Q^{(A)}$ with the generators in the q-commutator swapped.\\
\[ IQ^{(13)} = \frac{1}{q-q^{-1}} [Q^{(23)}, Q^{(12)}]_q -Q^{(1)} Q^{(3)} -Q^{(2)}Q^{(123)}\]
\[ IQ^{(24)} = \frac{1}{q-q^{-1}}[Q^{(34)}, Q^{(23)}]_q -Q^{(2)} Q^{(4)} -Q^{(3)}, Q^{(234)}\]
\[ IQ^{(124)} = \frac{1}{q-q^{-1}} [Q^{(123)}, Q^{(34)}]_q -Q^{(12)} Q^{(4)} -Q^{(3)} Q^{(1234)}\]
\[ IQ^{(14)} = \frac{1}{q-q^{-1}} [Q^{(234)}, Q^{(123)}]_q -Q^{(1)} Q^{(4)} -Q^{(23)} Q^{(1234)}\]
\[ IQ^{(134)} = \frac{1}{q-q^{-1}} [Q^{(12)}, Q^{(234)}]_q -Q^{(1)} Q^{(34)} -Q^{(2)}Q^{(1234)}.\]
We will refer to the 5 original generators as 'bosonic' and these 5 derived operators, and their involutions as 'fermionic.'

We will now demonstrate the algebra relations satisfied by these operators. First we discuss the case of commuting operators. 
\begin{prop}
	For distinct subsets $A,B \subset \{1, 2, 3, 4\}$, the operators $Q^{(A)}$ and $IQ^{(B)}$ commute whenever $A\subset B$, $B\subset A$, or $A\cap B=\emptyset$. 
\end{prop}
In the case that $A$ or $B$ are intervals,  i.e that the operators $Q^{(A)}$ or $Q^{(B)}$ are bosonic, the proof is simple using an earlier proposition. Otherwise, we use action on irreducible representations and the symbolic manipulator MAPLE to verify these relations. 

%
%
\begin{figure}
	\centering
	\includegraphics[scale=.4]{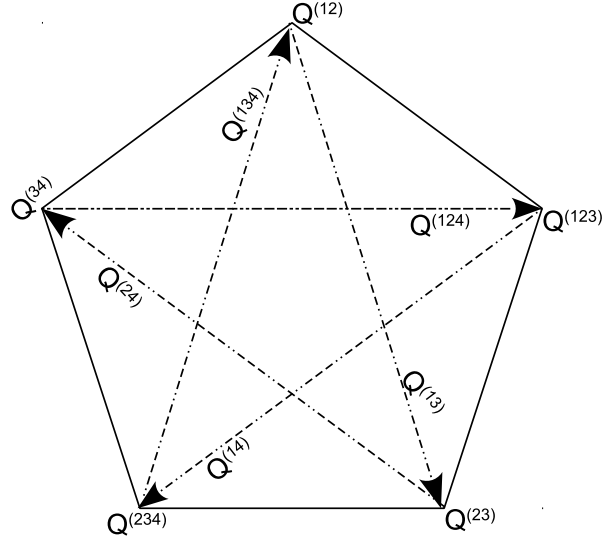}
	\caption{The symmetry compass. Vertices  represent the bosonic generators of the algebra. The dotted lines represent operators that do not commute and so give the remaining 10 linearly independent generators.  The orientation of the lines give the operators $Q^{(B)}$ and their involutions $IQ^{(B)}$ will have opposite orientations. }
	\label{AW4figuredirected}
\end{figure}
The remaining operations will come from generators that do not commute. From the previous theorem, two operators will not commute if their indices have non-trivial intersection and difference. We will now give the q-commutators of two such operators. However, since the q-commutator is not symmetric, we will need to define an ordering on the generators.   Let $i,j,k \subset \{1,2,3,4\}$ be nonempty and disjoint. So at most one of $i,j,k$ has two elements. Let there be equivalence of even permutations of $(i,j,k)$ (that is $(i,j,k)=(j,k,i)=(k,i,j)$). If all have only one element, then $(i,j,k)$ is an  allowable ordering if there is a permutation so that $i<j<k$. If one of these sets has 2 elements then up to a permutation $j$ has size 2, then $(i,j,k)$ is allowable if $i<k$.

This order structure can be seen on the symmetry compass as follows. For two operators, the ordering gives a q-commutator of the form $[Q^{A}, Q^{(B)}]_q$ if
\begin{itemize}
\item there is a directed edge from $Q^{(A)}$ to $Q^{(B)}$, in the case of two bosonic generators. 
\item the directed edge associated with $Q^{(B)}$ ends at $Q^{(A)}$ or the directed edge associated with $Q^{(A)}$ begins at $Q^{(B)}$, in the case that one is fermionic and the other bosonic. 
\item the two directed edges cross with the $Q^{(A)}$ ending clockwise of $Q^{(B)}$, in the case of two fermionic generators.  

\end{itemize}

\begin{prop}
Given an allowable ordering of the triple $(i,j,k)$, the operators $Q^{(ij)}$, $Q^{(jk)}$, and $Q^{(ik)}$ satisfy the symmetric $AW(3)$ relations,
\be \label{qcom}
\begin{aligned}	[]	
	\frac{1}{q-q^{-1}}[Q^{(ij)},Q^{(jk)}]_q &=Q^{(ki)} + I(Q^{(i)}Q^{(k)} + Q^{(ijk)}Q^{(j)})\\
	\frac{1}{q-q^{-1}}[Q^{(ki)},Q^{(ij)}]_q &=Q^{(jk)} + I(Q^{(k)}Q^{(j)} + Q^{(ijk)}Q^{(i)})\\
	\frac{1}{q-q^{-1}}[Q^{(jk)}, Q^{(ki)}]_q &=Q^{(ij)} + I(Q^{(j)}Q^{(i)} + Q^{(ijk)}Q^{(k)})\\
\end{aligned}
\ee
\end{prop}
Here we emphasize that $Q^{(ij)}=Q^{(ji)}.$ Thus, for each pair of distinct subsets $A,B \subset \{1,2,3,4\}$, either $Q^{(A)}$ commutes with $IQ^{(B)}$, or $i=A\setminus B$, $j=A\cap B$, and $k=B\setminus A$ are nonempty so $Q^{(A)}$, $Q^{(B)}$ are part of one of the above copies of $AW(3)$. In the case that the generators are all fermionic, the involution has a trivial action and we recover the $AW(3)$ algebra associated with the triple coupling $(ijk).$ Again, these algebra relations were verified using the action on irreducible representations via MAPLE. 

At this point we can say that the algebra $AW(4)$ closes in a fashion, described in the following theorem. 

\begin{prop}
The algebra $AW(4)$ contains 15 linearly independent generators: 5 bosonic generators $Q^{(A)}$, with $A$ a consecutive subset of $\{1,2,3,4\}$,  and $5$ fermionic generators,indexed by the remaining sets, along with their involutions.  The algebra relations are of two types: operators will either commute, or their q-commutators can be expressed as a linear combination of the generators, plus central terms. However, for the fermionic generators, the commutation relations will be between operators of different parities and the q- commutations relations between operators with the same parity. 
\end{prop}
This theorem is essentially a summary of the two previous, along with the observation that the involution on (\ref{qcom}) will result in q-commutator relations for the flipped generators.

\subsection{An additional identity involving triple q-commuators}
We finish this section, with a novel identity expressed solely in terms of the bosonic generators, i.e. those that arose naturally from the coupling schemes. These relations involve chosing triples of these bosonic generators from the symmetry compass,  $(A, B, C)$ which follow the ordering, i.e. go along two adjacent edges in the natural ordering or a single edge back and forth. 
\begin{figure}
	\centering
	\includegraphics[scale=.4]{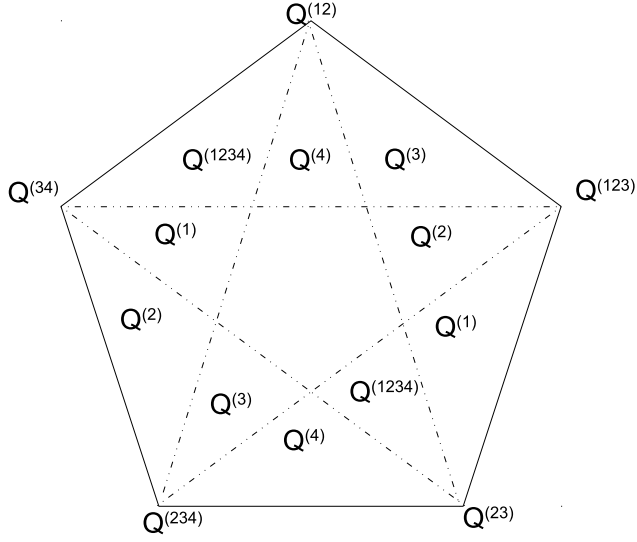}
	\caption{The symmetry compass, with central terms.   Two non-commuting operators generated a subalgebra isomorphic to $AW(3)$ with central elements including the third vertex of the triangle along with the 3 elements within that triangle.  }
	\label{AW4figure}
\end{figure}

\begin{prop}
The operators $Q^{(\mu)}$ arising from the tensor product of 4 irreducible representations of $U_q(su(1,1))$ satisfy 
\be \begin{split}
	[[Q^{(A)},Q^{(B)}]_q,Q^{(C)}]_q+[[Q^{(\alpha)}, Q^{(\beta)}]_q, Q^{(\gamma)}]_q+[[Q^{(X)}, Q^{(Y)}]_q, Q^{(Z)}]_q
	= \\ [[Q^{(A)}, Q^{(\beta)}]_q, Q^{(Z)}]_q+[[Q^{(X)}, Q^{(B)}]_q, Q^{(\gamma)}]_q+[[Q^{(\alpha)}, Q^{(Y)}]_q, Q^{(C)}]_q\, ,
\end{split} \label{AWMaster} \ee
with the choice of triples as given in Tables \ref{3nonc} and \ref{2nonc} below. 
\end{prop}
%

\begin{table}[h]
\centering
\begin{tabular}{c|c|c}
(A, B, C) & $(\alpha, \beta, \gamma)$ & (X, Y, Z)\\
\hline
(234, 12, 23) & (1, 2, 4 )& (0, 34,  123)\\
(23, 12, 234) & (2, 1, 4) & (0, 123, 34)\\
(34, 123, 234) & ( 4, 1234, 2 )& (0, 12,  23)\\
(234, 123, 34) & ( 1234, 4, 2 )& (0, 23, 12)\\
(12, 23, 34) & (2, 3, 1234) & (0, 123, 234)\\
(34, 23, 12) & (3, 2, 1234) & (0, 234, 123)\\
(123, 234, 12)& (1234, 1, 3 )& (0,23, 34)\\
(12, 234, 123)& (1, 1234, 3)& (0,34, 23)\\
(23, 34, 123)& (3, 4, 1 )& (0,234, 12)\\
(123, 34, 23)& (4, 3, 1 )& (0,12, 234)
\end{tabular}
\caption{Choices of parameters for 3 noncommuting operators. }
\label{3nonc}
\end{table}
\begin{table}[h]
\centering
\begin{tabular}{c|c|c}
(A, B, C) & $(\alpha, \beta, \gamma)$ & (X, Y, Z)\\
\hline
(12, 23, 12) & (3, 2, 0) & (0, 1, 123)\\
(23, 12, 23)& (1, 2, 0) & (0, 3, 123)\\
(123, 234, 123) & (1, 1234, 0) & (0, 4, 23)\\
(234, 123, 234) & (4, 1234, 0) & (0, 1, 23)\\
(23, 34, 23) & (4, 3, 0) & (0, 2, 234)\\
(34, 23, 34) & (2, 3, 0) & (0, 4, 234)\\
(234, 12, 234) & (2, 1, 0) & (0, 1234, 34)\\
(12, 234, 12) & ( 1234, 1, 0) &(0, 2, 34)\\
(34, 123, 34)& (1234, 4, 0) & (0, 3, 12)\\
(123, 34, 123) & ( 3, 4, 0) & (0, 1234, 12) 
\end{tabular}
\caption{Choices of parameters for triple commutators of pairs of operators. }
\label{2nonc}
\end{table}
\noindent This proposition is proven by direct computation of the relations using the discrete series representation discussed above. The central  operator $Q^{(0)}\equiv -I$ is introduced for convenience. 

To read off the choices of coefficients in Table \ref{3nonc} from  Figure \ref{AW4figure}, we note that each triple is associated with triangle determined by a pair of dotted lines and a solid line. The choice of $(\alpha, \beta, \gamma)$ is determined by following the lines along the triangle and picking the center coefficients from the subalgebras generated by the lines parallel to the sides of the triangle. The choice of $Y $ and $Z$ is just the edges of the dotted line parallel to the solid line of the triangle ordered following the orientation of the triangle.  For the commutators of two operators, given in Table \ref{2nonc}, you read off the coefficients $(\alpha, \beta, 0)$ and $(0, Y, Z)$ by going around the triangle for the copy of $AW(3)$ beginning with the corner $Q^{(A)}$.

\section{Conclusion}
In summary, we have use the coproduct structure of the universal enveloping algebra $U_q(sl(2,\mathbb{C})$ to define generators for the well-known quadratic algebras $AW(3)$ and its novel extension $AW(4)$ and give their structure equations.  The algebra $AW(4)$ is comprised of several copies of $AW(3)$, generated by q-commutators along with several commutator relations. Unlike the $AW(3)$ case, there is a non-trivial inversion leading to a pairity structure on the algebra that must be respected by the algebra relations.

There are many areas where this research can be extended. First, the  extension of the generators to arbitrary dimension is clear and  will be the subject of further investigation.  By construction, these algebras as associated with coupling schemes of irreducible representations of the underlying quantum algebra and in certain cases the interbasis expansion coefficients are known to be given by q-Racah polynomials\cite{genest2017coupling}.   Another subject of further investigation is into other decompositions and their associated orthogonal functions. Finally, we mention that the symmetry inherent in the construction, including permutations of indices and duality in the interbasis expansion, is a subject that has only begun to be explored and mined for special function identities. \\

\noindent{\bf Acknowledgements:} SP acknowledges  The Simons Foundation Collaboration grant \# 3192112 for support of this research.

\end{document}